\documentclass[a4paper,twoside,12pt]{amsart}
\usepackage{amsfonts,amssymb}
\newcommand\pic{\mathop{\rm Pic}\nolimits}
\newcommand\sing{\mathop{\rm Sing}\nolimits}
\renewcommand{\P}{{\mathbb P}}

\begin{document}

\author{Sonia Brivio}
\address{S. Brivio: Dipartimento di Matematica,
Universit\'a di Pavia, via Ferrata 1, 27100 Pavia}
\email{sonia.brivio@unipv.it}
\author{Alessandro Verra}
\address{A. Verra: Dipartimento di Matematica,  Terza Universit\'a di Roma,
Largo S.Murialdo, 00146 Roma} \email{verra@mat.uniroma3.it}
\date{}

\title{The Brill-Noether curve of a stable vector bundle on a genus two curve.}
 
\footnotetext[1]{\noindent Partially supported by the research programs {\em
'Moduli Spaces and Lie Theory'} and {\em 'Geometry on  Algebraic  Varieties'}
of the Italian Ministry of Education.}
\footnotetext[2]{\noindent 2000 Mathematics Subject Classification: 14H60}
\begin{abstract} Let ${\mathcal U}_r$ be the moduli space of rank $r$ vector bundles with trivial determinant on a smooth 
curve of genus $2$.  The map 
${\theta}_r \colon {\mathcal U}_r \to \vert r {\Theta} \vert$, which associates to a general bundle its theta divisor, is generically finite. In this paper we give a geometric interpretation 
of the generic fibre of ${\theta}_r$.
\end{abstract}
\newtheorem{teo}{Theorem}[section]
\newtheorem{prop}[teo]{Proposition}
\newtheorem{lem}[teo]{Lemma}
\newtheorem{rem}[teo]{Remark}
\newtheorem{defi}[teo]{Definition}
\newtheorem{ex}[teo]{Example}
\newtheorem{cor}[teo]{Corollary}
 
\maketitle 

\section{Introduction. } In this note we deal with the moduli space ${\mathcal 
U}_r$ of semistable vector bundle of rank $r$ and degree $r(g-1)$ over  a
smooth, irreducible complex projective curve of genus $g \geq 2$. ${\mathcal  U}_r$
is endowed with the Brill-Noether locus
\[ {\Theta}_r:= \{ \ [E] \in {\mathcal  U}_r \mid h^0(E) \geq 1 \ \}
\] which is an integral Cartier divisor and it is known as the generalized theta
divisor of ${\mathcal  U}_r$, see \cite{Drezet}, \cite{BNR}. Moreover the tensor
product defines a morphism
\[ f: {\mathcal  U}_r \times \pic^0(C) \to {\mathcal  U}_r
\] and we can consider the pull-back $f^*\Theta_r$ of $\Theta_r$. Let
$[E] \in {\mathcal  U}_r$ be the moduli point of the vector bundle
$E$ and let $\det\ E \cong M^{\otimes r}$, it is well known that then
\[ {\mathcal  O}_{[E] \times \pic^0(C)}(f^* \Theta_r) \cong {\mathcal 
O}_{\pic^0(C)}(r\Theta_M)
\]   where 
\[ 
{\Theta}_M := \{ N \in \pic^0(C) \mid h^0(M \otimes N) \geq 1 \}.
\]
 Note that $M$ is a line bundle of degree $g-1$ and that
$\Theta_M$ is a theta divisor on $\pic^0(C)$. We define
\begin{equation*} 
{\Theta}_{E} := f^*\Theta_r \cdot [E] \times \pic^0(C) 
\end{equation*} 
if the intersection is proper. In this case we will say that \it
$\Theta_E$ is the theta divisor of $E$. \rm The construction of
$\Theta_E$ allows us to define a rational map as follows. Consider 
\begin{equation*} 
{\mathcal  T}_r \ := \ \bigcup_{M \in \pic^{g-1}(C)}  |
r\Theta_M |, 
\end{equation*}  
it is a  standard fact that ${\mathcal  T}_r$ has a natural
structure of projective bundle over $\pic^{r(g-1)}(C)$. So we omit its
construction, we only mention that the corresponding projection 
\[
p: {\mathcal  T}_r \to \pic^{r(g-1)}(C)
\]
is defined as follows: $p(D) = M^{\otimes r}$ iff $D \in \ |  r\Theta_M
| $. Notice that the only elements of multiplicity $r$ in
$| r\Theta_M |$ are exactly the divisors $r\Theta_{M \otimes
\eta}$, where $\eta$ varies in the set of the elements of order $r$ in
$\pic^0(C)$. Therefore the map $p$ is well defined. In the following we will
study the rational map
\begin{equation*} 
{\theta}_r: {\mathcal  U}_r \to {\mathcal  T}_r 
\end{equation*} 
which
associates to a general $[E] \in {\mathcal  U}_r$ the corresponding theta divisor
${\Theta}_E \in {\mathcal  T}_r$. Let
\[ 
\det: {\mathcal  U}_r \to \pic^{r(g-1)}(C)
\] 
be the  determinant map, it is well known that  ${\mathcal  T}_r$ is the
projectivization of $ \det_* {\mathcal  O}_{{\mathcal  U}_r}({\Theta}_r)^*$ and that
$\theta_r$ is the induced tautological map. In particular it follows that $p
\cdot \theta_r = \det$. \hfill \par \noindent We will say that $\theta_r$ is \it
the theta map. \rm \par \noindent Too many questions are still unsettled about
the theta map, excepted for the case $r \leq 2$: see e.g. \cite{Beauville2} for
a general survey. This situation is probably related to the fact that the next
basic question is still mostly unsolved.
\hfill \par \hfill \par \noindent 
\bf QUESTION \rm {\it Is $\theta_r$ generically finite onto its image?}
\rm \bigskip \par \noindent Actually the main difficulty here is that
$\theta_r$ is not a morphism in most of the cases \cite{Raynaud}. Thus, in
spite of the ampleness of $\Theta_r$, it is not a priori granted that
$\theta_r$ is generically finite onto its image. \hfill \par \noindent In this
paper we give a natural geometric interpretation of the fibres of the map
$\theta_r$ for a curve $C$ of genus two. A very special feature of this case is
that
\[ 
\dim \ {\mathcal  U}_r = \dim {\mathcal  T}_r = r^2 + 1,
\]
so the generic finiteness of $\theta_r$ is even more expected. Applying our
description of the fibres we prove the generic finiteness of
$\theta_r$. \par \noindent Such a result is  not new: Beauville recently proved
it using a different, relatively simple method, see
\cite{Beauville1}. We believe that our description has some interest in itself
and we hope to  use it for further applications, in particular to compute the
degree of $\theta_r$. 
\hfill \par \noindent Our approach relies on Brill-Noether theory for curves
contained in a genus two Jacobian. Let $D \in \theta_r({\mathcal  U}_r)$ be a
sufficiently general element, then $D$ is a smooth curve of genus $r^2+1$ in
$\pic^0(C)$: see section 2. Consider the Brill-Noether locus
\[
 W^{r-1}_{r^2}(D) \ = \ \{ L \in \pic^{r^2}(D) \mid h^0(L) \geq r \}
\]
and observe that its expected dimension is one, in other words the
Brill-Noether number $\rho(r-1,r^2,r^2+1)$ is one.  Our main result can be
summarized as follows:
\hfill\par \hfill\par \noindent
\bf THEOREM \rm  \it {\it Each point $[E] \in \theta_r^{-1}(D)$ defines an
irreducible component
\[
 C_E \subset W^{r-1}_{r^2}(D)
\]
biregular to $C$. Let $Z$ be the set of all irreducible components of
$W^{r-1}_{r^2}(D)$ and let
\[
 i_D: \theta_r^{-1}(D) \to Z
\]
be the map sending $[E]$ to $C_E$. Then $i_D$ is injective. \rm
\hfill\par \hfill\par \noindent The statement clearly implies that $\theta_r$
is generically finite. We define $C_E$ as the \it Brill-Noether curve of $E$.
\rm Fixing appropriately a Poincar\'e bundle ${\mathcal  P}$ on $D \times
\pic^{r^2}(D)$ it turns out that $E$ is the restriction of
$\nu_* {\mathcal  P}$ to $C_E$, where $\nu$ is the projection onto
$\pic^{r^2}(D)$. In particular the family of the fibres of $E$ is just the
family of the spaces
$H^0(L)$, $L \in C_E$. Notice that the choice of $\mathcal P$, hence of $\det \
E$, depends on the embedding $D \subset \pic^0(C)$ and it is essentially
explained in the final part of this note. \par \noindent  To have a typical
example of what happens, the reader can consider the case $r = 2$. In this case
$D$ is a curve of genus 5  endowed with a fixed point free involution which is
induced by the $-1$ multiplication of $\pic^0(C)$.  Since $r = 2$ the
Brill-Noether locus
$W^1_4(D)$ is exactly the singular locus of the theta divisor of $\pic^4(D)$. It
follows from the theory of Prym varieties that $W^1_4(D)$ is the union of two
irreducible curves: one of them has genus 4, the other one is just a copy of
$C$,(see also \cite{Teix}). This is the Brill-Noether curve of a stable rank
two vector bundle $E$ such that $\theta_2([E]) = D$. In higher rank the general
theory of Prym-Tjurin varieties can certainly provide further information on
$W^{r-1}_{r^2}(D)$ and hence on the fibres of $\theta_r$. However, in order to
get them, a very explicit description is needed for the Prym-Tjurin
realizations of a genus two Jacobian. \par \noindent  On Jacobians of higher
genus several extensions of the above constructions are  possible and perhaps
deserve to be considered in the study of the theta maps. We hope to have
underlined with this note the multiplicity of the links between moduli of
vector bundles on a curve $C$, Prym-Tjurin realizations of its Jacobian $JC$
and Brill-Noether theory for curves in $JC$.
\hfill\par We wish    to thank the referee for some helpful comments.
\bigskip \par \noindent \it The second author wishes to express his gratitude
to Jacob Murre, on the occasion of the Proceedings of the Conference in his
honour, for the friendship and for the scientific attention received during
three decades.
\bigskip \noindent \rm
\section{Notations and preliminary results.}  From now on $C$ is a smooth,
irreducible, complex projective curve of genus $2$. Let $C^{(2)}$ be the
2-symmetric product of $C$, a point of such a surface is a divisor $x+y$ with
$x,y \in C$. We consider the map
\begin{equation*} 
a: C^{(2)} \to \pic^0(C)
\end{equation*} 
sending $x+y \in C^{(2)}$ to $\omega_C(-x-y)$. Of course $a$ is
the composition of the Abel map defined by $\omega_C$ with
$-1$ multiplication on $\pic^0(C)$. Therefore $a = -\sigma$, where $\sigma:
C^{(2)} \to \pic^0(C)$ is the blowing up of the zero point. For each fibre $|
r\Theta_M |$ of the projective bundle ${\mathcal  T}_r$ we have the  linear
isomorphism
\[ a^*_M: | r\Theta_M | \to | a^*r\Theta_M |
\]
defined by the pull-back. Let $\Theta_E$ be the theta divisor of $[E]$, we
will keep the following notation 
\begin{equation*} D_E := a^* \Theta_E.
\end{equation*} 
$D_E$ is an effective divisor in $C^{(2)}$ which is supported on
the set
\begin{equation*}
\lbrace x+y \in C^{(2)} \mid h^0(E \otimes \omega_C(-x-y)) \geq 1 \rbrace.
\end{equation*}
$D_E$ is biregular to $\Theta_E$ if the zero point is not in $\Theta_E$,
otherwise $D_E$ is the union of the projective line 
$| \omega_C |$ and of a curve  birational to $\Theta_E$. For our
convenience we are more interested to $D_E$ than to $\Theta_E$. 
\par \noindent At first we want to prove that a general $D_E$ is smooth, to do
this we need Laszlo's singularity theorem, see
\cite{Laszlo}:
\hfill \par \hfill \par
\begin{teo} \label{twoone} The multiplicity of $\Theta_r$ at its stable point
$[E]$ is $h^0(E)$.
\end{teo}  
\hfill \par \par
\begin{prop}  \label{twotwo} Let $[E]$ be a general stable point of ${\mathcal 
U}_r$ then  
\[ h^0(E \otimes \omega_C(-x-y)) \leq 1, \ \ \forall x+y \ \in \
C^{(2)}.\] 
\end{prop}
\begin{proof} By induction on $r$. Let $r = 1$ then $D_E = C$ and $E$ is a
general line bundle of degree 1, in particular $| E \otimes
\omega_C |$ is a base-point-free pencil and this implies the statement. Let
$r \geq 2$, we can assume by induction that there exist general $[B] \in
{\mathcal  U}_{r-1}$ and  $[A] \in {\mathcal  U}_1$ $=$ $\pic^1(C)$ satisfying
the statement. We consider the exact sequence
\[
 0 \to B \to E \to A \to 0
\]
defined by the vector $e \in \mbox{\rm \mbox{\rm \mbox{\rm \mbox{\rm
Ext}}}}^1(A,B)$. Tensoring such a sequence by $\omega_C(-x-y)$ and passing to
the long exact sequence we obtain the coboundary map
\[ 
e_{x+y}: H^0(\omega_C \otimes A(-x-y)) \to H^1(\omega_C \otimes B(-x-y)).
\]
{\bf Claim } The statement holds for $E$ iff $e_{x+y}$ has maximal rank for
every $x+y$. \hfill\par \noindent Proof: We have $h^0(\omega_C \otimes A(-x-y))
\leq 1$ and 
\[
h^0(\omega_C \otimes B(-x-y))=h^1(\omega_C \otimes B(-x-y))\leq 1.
\]
 Then the statement follows from the above mentioned long exact
sequence. \medskip \par \noindent Finally it is obvious that $e_{x+y}$ has
maximal rank except possibly for points $x+y$ with  $h^0(\omega_C \otimes
A(-x-y))= h^0(\omega_C
\otimes B(-x-y))= 1$. The set of these points is $D_A \cap D_B$. Since $A$
is general we can assume that $D_A \cap D_B$ is finite. Let
$x+y \in D_A \cap D_B$ then $e_{x+y}$ has not maximal rank iff it is the zero
map.  It is  a standard property that, in the present case, the locus \rm
\[ H_{x+y} = \{ e \in \mbox{\rm \mbox{\rm \mbox{\rm \mbox{\rm \mbox{\rm
Ext}}}}}^1(A,B) \
\mid \ e_{x+y} \ \text {is the zero map} \}
\] 
is a hyperplane. Let $H := \bigcup H_{x+y}$, $x+y \in D_A \cap D_B$, then a
general $e \in \mbox{\rm \mbox{\rm \mbox{\rm \mbox{\rm Ext}}}}^1(A,B) - H$
defines a semistable $E$ satisfying the condition of the statement. Since this
condition is open on ${\mathcal  U}_r$ the result follows. 
\end{proof}
\begin{cor} \label{twothree} Let $[E]$ be a general stable point of ${\mathcal 
U}_r$, then $D_E$ is smooth.
\end{cor}
\begin{proof} Let $f: {\mathcal  U}_r \times \pic^0(C) \to {\mathcal  U}_r$ be
the map defined via tensor product, recall that  
\[ \Theta_E = f^* \Theta_r \cdot [E] \times \pic^0(C) \subset {\mathcal  U}_r
\times \pic^0(C). \] 
Therefore $\Theta_E$ is the fibre of the projection $q: f^*\Theta_r \to
{\mathcal  U}_r$. Then, by generic smoothness, a general $\Theta_E$ is smooth
if $\Theta_E \cap \sing \ f^*\Theta_r = \emptyset$. On the other hand
$f$ is smooth, with fibres biregular to $\pic^0(C)$. The smoothness of $f$
implies that $\sing \ f^*\Theta_r$ $=$  $f^* \sing \ \Theta_r$. Therefore, by
Laszlo's singularity theorem and the definition of $f$, we have
\[\sing \ f^*\Theta_r = \lbrace ([E],\xi) \in {\mathcal  U}_r \times \pic^0(C)
\mid h^0(E(\xi)) \geq 2 \rbrace.
\] 
But the previous proposition implies that $h^0(E(\xi)) \leq 1$, for all  
$\xi \in \pic^0(C)$. Then it follows that $\Theta_E \cap \sing \ f^*
\Theta_r = \emptyset$ and hence a general $\Theta_E$ is smooth. The same holds
for $D_E$. 
\end{proof}

\section{The tautological model $P_E$} 
Now we want to see that the above curve
$D_E$ appears as the singular locus of some natural tautological model of
$\P E^*$ in $\P^{2r-1}$. 
\begin{prop}  \label{threeone} Let $E$ be any semistable point of ${\mathcal 
U}_r$ then  
 \par\noindent {\rm 1)} $h^1(\omega_C \otimes E) = 0$ and $h^0(\omega_C \otimes
E) = 2r$.\par
\noindent
\par\noindent {\rm 2)}  $\omega_C \otimes E$ is globally generated unless $E$
is not stable and
$\mbox{\rm Hom}(E,{\mathcal  O}_C(x))$ is non zero for some point $x
\in C$. 
\end{prop}
\begin{proof}  1)  \rm  By Serre duality $h^1(\omega_C \otimes E) =
h^0(E^*)$. Since $E^*$ is semistable of slope $-1$ it follows
$h^0(E^*) = 0$. Then we have $h^0(\omega_C \otimes E) = 2r$ by Riemann-Roch.
\hfill \par \noindent
  2)  \rm By   1) \rm $E$ is globally generated iff $h^0(\omega_C \otimes
E(-x)) = r$, $\forall \  x \in C$. By Serre duality this is equivalent to
$\mbox{\rm Hom}(E, {\mathcal  O}_C(x)) = 0$, $\forall \ x \in C$. Notice also
that 
$\mbox{\rm Hom}(E,{\mathcal  O}_C(x)) \neq 0$ implies that $E$ is not stable.
This completes the proof. \end{proof} 
In this section we assume that $[E] \in {\mathcal  U}_r$ has the
following properties (satisfied by a general $[E]$):  \par
\noindent \it  - $\omega_C \otimes E$ is  globally generated,  \hfill \par
\noindent - $D_E$ exists i.e. $[E]$ is not in the indeterminacy locus of
$\theta_r: {\mathcal  U}_r \to {\mathcal  T}_r$, \hfill \par \noindent - $D_E$ is smooth,
\rm \hfill \par \noindent To
simplify our notations we put
\begin{equation*} F := \omega_C \otimes E \ \ \text {and} \ \ \P_E :=
\P F^*.
\end{equation*}
\lem \label{threetwo}
Let $F$ be general and let $\overline F$ be defined by the
standard exact sequence
\[
 0 \to F^* \to H^0(F)^* \otimes {\mathcal  O}_C \to \overline F \to 0
\] 
induced by the evaluation map. Then $\overline F$ is stable. In particular
the map
\[ j: {\mathcal  U}_r \to {\mathcal  U}_r
\] 
sending $[\omega_C ^{-1} \otimes F]$ to $[\omega_C ^{-1} \otimes \overline
F]$ is a birational involution. 
\endlem
\begin{proof} First of all we claim that $\overline F$ is semistable for $F$ general enough. 
Let $F_o = L^r$, where  $L \in \pic^3(C)$ is globally generated,
then $h^0(F_o) = 2r$ and $\overline {F}_o = F_o$. Up to a base change there
exists an integral variety $T$ and a vector bundle $\mathcal F$ over $T \times
C$ such that the family of vector bundles
$\lbrace F_t := {\mathcal  F} \otimes {\mathcal  O}_{t \times C}, \ t \in T \rbrace$
dominates ${\mathcal  U}_r$ and contains $F_o$.  By semicontinuity we can assume, up
to replacing $T$ by a non empty open subset, that $h^0(F_t) = h^0(F_o)= 2r$ and $F_t$ is globally generated.
So it is standard to construct  from ${\mathcal  F}$ a vector bundle
$\overline{\mathcal  F}$ on $T \times C$  with the following property: 
$\overline{\mathcal  F} \otimes {\mathcal  O}_{t \times C}$ $=$
$\overline {F}_t$, for each $t \in T$. Since $\overline {\mathcal  F}  \otimes
{\mathcal  O}_{o\times C}$ $=$ $L^r$ is semistable, the same holds for a
general vector bundle $\overline {\mathcal  F} \otimes {\mathcal  O}_{t \times
C}$. Hence the claim follows. 
Let $F$ be a general stable  bundle: $\overline F$ is semistable,  by lemma (3.1)
$h^0(\overline F) =2r$, moreover since $h^0(F^*)= 0$, we have  $H^0(F)^* \simeq H^0(\overline F)$ and 
$\overline F$ is globally generated. So $j$ is defined at $\overline F$, actually $j(\overline F)= F$. 
This implies that $j$ is a birational involution and $\overline F$ is stable too. 
\end{proof} 
\rm Since $F$ is globally generated the map defined by ${\mathcal 
O}_{\bold P_E}(1)$ is a morphism
\begin{equation*} 
u_E: \P_E \to \P^{2r-1} := \P H^0(F)^*.
\end{equation*} 
In particular the restriction of $u_E$ to any fibre $\P_{E,x}$ of $\P_E$ is a
linear embedding
\[ 
u_{E,x}: \P_{E,x} \to \P^{2r-1}.
\]
\begin {defi} \rm The image of $u_E$, (of $u_{E,x}$), will be denoted 
 $P_E$,($P_{E,x}$). \end{defi}  
\noindent For any  $ \ d \in C^{(2)}$,   $F_d:= F \otimes {\mathcal  O}_d$ can be
naturally seen as a rank $r$ vector bundle over $d$. Note that its
projectivization is $p^*d$, where  $p: \P_E \to C$ is  the projection
map. In particular the evaluation map $e_d: H^0(F_d) \otimes {\mathcal  O}_d
\to F_d$ defines an  embedding 
\[ i_d: p^*d \to \P H^0(F_d)^*.
\] 
We have $\P H^0(F_d)^* = \P^{2r-1}$ and moreover $i_d(p^*d)$ is
the union of two disjoint linear spaces of dimension $r-1$ if $d$ is  smooth.
The next lemma is therefore elementary.
\begin{lem} \label{threefour} Let $o \in \P H^0(F_d)^*$ be a point not in
$i_d(p^*d)$, then there exists exactly one line
$L$ containing $o$ and such that $Z := i_d^*L$ is a 0-dimensional scheme of
length two. Moreover let $\lambda$ be the linear projection of centre $o$, then
$Z$ is the unique 0-dimensional scheme of length two on which $\lambda \cdot
i_d$ is not an embedding. 
\end{lem} \noindent The central arrow in the long exact sequence 
\begin{equation} \label{nine}
0 \to H^0(F(-d)) \to H^0(F) \to H^0(F_d) \to
H^1(F(-d))
\to 0 
\end{equation}  defines a linear map 
\[
\lambda_d: \P H^0(F_d)^* \to \P H^0(F)^* = \P^{2r-1}, 
\]  
and from the construction clearly $\lambda_d \cdot i_d$ is
the map 
\[ u_E|_{p^*d}: p^*d \to \P^{2r-1}.
\]
\begin{prop} \label{threefive} $u_E|_{p^*d}$ is not an embedding if and only if
$d
\in D_E$. 
\end{prop}
\begin{proof} From the previous remarks and lemma  \ref{threefour} it follows
that $u_E|_{p^*d}$ is not an embedding iff $\lambda_d$ is not an isomorphism. By
the long exact sequence (\ref{nine}) $\lambda_d$ is not an isomorphism iff
$h^0(E(-d))
\geq 1$, that is if $d \in D_E$. \end{proof} 

We want to use the previous
results to study the singular locus of $P_E$. Let $\mbox{\rm Hilb}_2(\P_E)$ be
the Hilbert scheme of 0-dimensional schemes $Z \subset \P_E$ of length two, we
simply consider  its closed subset
\begin{equation*}
\Delta = \lbrace Z \in \mbox{\rm Hilb}_2(\P_E) \mid u_E|_Z \ \text {is not an
embedding} \rbrace.
\end{equation*}     
Let $ Z \in \Delta$ then $Z \subset p^*d$, where $d := p_*Z$
belongs to $D_E$. So we have a morphism
\begin{equation*} 
p_*: \Delta \to D_E
\end{equation*} 
sending $Z$ to $d$. 
\begin{prop} \label{threesix} Let $E$ be general then $p_*: \Delta \to D_E$ is
biregular.
\end{prop}
\begin{proof} Let $Z \in \Delta$ and let $p_*Z = d$, then $Z$ is embedded in
$p^*d$. Since $d \in D_E$ we have $h^0(F(-d)) = 1$, see prop. \ref{twotwo}. This
implies that the linear map
\[
\lambda_d: \P H^0(F_d)^* \to \P^{2r-1}
\] 
is the projection of centre a point $o$ with image a hyperplane in
$\P^{2r-1}$. Then, by lemma \ref{threefour}, $Z$ is the unique element of
$\Delta$ which is contained in $p^*d$. Hence $p_*$ is injective. Conversely let
$d \in D_E$, then $\lambda_d \cdot i_d$ is not an embedding on exactly one
0-dimensional scheme $Z \subset p^*d$ of length two. Since $\lambda_d \cdot i_d
= u_E|_{p^*d}$, it follows that $Z$ is in $\Delta$ and that $p_*$ is surjective.
Since $D_E$ is a smooth curve, $p_*$ is biregular.
\end{proof}
\begin{prop} \label{threeseven} Assume $r \geq 2$ and $E$ general, then  $u_E:
\P_E  \to P_E$ is the normalization map and $\sing \ P_E$ is an irreducible
curve.
\end{prop}
\begin{proof} Let $\tilde D \subset \P_E$ be the image of the curve  
\[
\tilde \Delta = \lbrace (Z,q) \in \Delta \times \P_E \mid q \in \mbox{\rm
Supp} \ Z \rbrace \]
under the projection $\Delta \times \P_E \to \P_E$. The set $\tilde D$ is the
locus of points where $u_E$ is not an embedding: it is a proper closed set as
soon as $r \geq 2$. Hence $u_E: \P_E \to P_E$ is a morphism of degree one if $r
\geq 2$. Since $\P_E$ is smooth, $u_E$ is the normalization map if each of its
fibres is finite. Assume $u_E$ contracts an irreducible curve $B$ to a point
$o$. $B$ cannot be in a fibre $\P_{E,p}$: otherwise $D_E$ would contain the
curve $\lbrace z+p, z \in C \rbrace$ and would be reducible. Hence $o \in \cap
P_{E,x}, \ x \in C$. Let $x+y \in C^{(2)}$ with $x \neq y$, then $P_{E,x} \cup
P_{E,y}$ is contained in a hyperplane and hence $h^0(E
\otimes \omega_C(-x-y)) \geq 1$. This implies $D_E = C^{(2)}$: a contradiction.
It remains to show that $\sing  P_E$ is an irreducible curve: this is clear
because $\sing \ P_E = u_E(\tilde D)$ 
\end{proof}

\section {The line bundle $H_E$} We will keep the generality assumptions and
the notations of the previous section. Recall that $d \in D_E$ uniquely defines
a 0-dimensional scheme $Z_d \subset p^*d$ of length two such that $u_E|_{Z_d}$
is not an embedding, in particular $u_E(Z_d)$ is a point.
\defi $h_E: D_E \to \P^{2r-1}$ is the morphism sending $d$ to
$u_E(Z_d)$, moreover
\[
 H_E := h_E^* {\mathcal  O}_{\P^{2r-1}}(1).
\] 
\enddefi \rm \noindent
\rem \rm \it 1. \rm Let $F := \omega_C \otimes E$ and let $q_1, q_2 : C \times
C \to C$ be the projections. Note that
$q_1^*F \oplus q_2^*F$ descends to a vector bundle $F^{(2)}$ on $C^{(2)}$ via
the quotient map $C \times C \to C^{(2)}$.  Moreover the evaluation $H^0(F) \to
F_x \oplus F_y$ induces a natural map
\[ 
e: H^0(F) \otimes {\mathcal  O}_{C^{(2)}} \to F^{(2)}.
\] 
Then $D_E$ is the degeneracy locus of $e$ and $H_E$ is its cokernel. This
implies that the sheaf $H_E$ can be defined for  every  curve $D_E$ and that
$H_E$ is a line bundle iff $h^0(\omega_C \otimes E(-x-y)) = 1$ for each $x+y
\in D_E$.

\par \noindent \it 2. \rm A very simple geometric definition  of $h_E$ can be
given as follows: let $d = x+y \in D_E$ with $x \neq y$ then 
\begin{equation*} 
h_E(d) = P_{E,x} \cap P_{E,y} = u_E(Z_d).
\end{equation*}
\endrem
\begin{prop} \label{fourthree} $h_E: D_E \to \P^{2r-1}$ is generically injective
if $E$ is general and $r \geq 2$.
\end{prop}
\begin {proof} Let $  U = \lbrace x+y  \in D_E \mid x \neq y
\rbrace$, assume that $d_1, d_2 \in U$ and  $h_E(d_1) = h_E(d_2) = o$: then  
$ o \in \bigcap_{i = 1 \dots 4} P_{E,x_i}$, where $\Sigma \  x_i = d_1 + d_2$. 
We consider the standard exact sequence
\[
 0 \to F^* \to H^0(F)^* \otimes {\mathcal  O}_C \to \overline F \to 0,
\]
where $F = \omega_C \otimes E$. The long exact sequence identifies
$H^0(F)^*$ to a subspace of $H^0(\overline F)$. Hence $o$ is a 1-dimensional
space generated by some $s \in H^0(\overline F)$. It is standard to verify that
$o \in \bigcap_{i = 1 \dots 4} P_{E_{x_i}}$ iff $s$ is zero on $d_1+d_2-d$,
where $d$ is the M.C.D. of $d_1,d_2$. By Lemma~\ref{threetwo} $\overline F$ is
stable, hence we must have $deg \ d \geq 2$ that is $d_1 = d_2 = d$. 
\end{proof}
 \begin{prop} \label{fourfour} $H_E$ has degree $r^2 + 2r$. \end{prop}
\begin{proof} Set theoretically we have $h_E(D_E) = \sing \ P_E$, hence $\sing \
P_E$ is an irreducible curve. Let $o = h_E(x+y)$ be general then $x \neq y$ and
moreover $u_E^*(o)$ is supported on two closed points $o'$ and $o''$: this
follows because $h_E$ is generically  injective. 
\par {\bf Claim}: \emph{
The tangent map $du_E$ is injective at $o'$ and $o''$.}
\par
Let $\widetilde{D(u_E)}$ be the double point scheme of
$u_E$, defined as in \cite[ p.~166]{Fulton}. $\widetilde{D(u_E)}$ is contained
in
$\widetilde{\P_E \times \P_E}$, where $\pi:
\widetilde{\P_E \times \P_E} \to \P_E \times \P_E$ is the blowing up of the
diagonal $\Delta$. A point of $\widetilde{D(u_E)}$ is either the inverse image
by $\pi$ of a pair $(o',o'')$ in $\P_E \times \P_E - \Delta$ such that $u_E(o')$
$=$ $u_E(o'')$ or it is a point in $\pi^{-1}(\Delta)$ parametrizing a 1
dimensional space of tangent vectors to $\P_E$ on which $du_E$ is zero. In the
former case we have also $p(o') \neq p(o'')$ because $u_E$ is injective on each
fibre of $p$. On the other hand it is clear that, in our situation,
\[
 q^{-1}(D_E) = (p \times p) \cdot \pi (\widetilde{D(u_E)})
\]
where $q: C \times C \to C^{(2)}$ is the quotient map. Thus $du_E$ is not
injective at most along fibres $\P_{E,z}$ such that $2z \in D_E$. $D_E$
cannot be the diagonal of $C^{(2)}$ because $D_E^2 = 2r^2$. Hence $D_E$
contains finitely many points
$2z$ and we can choose the above point $o = h_E(x+y)$ so that $2x$ and $2y$ are
not in $D_E$. This implies our claim.
\par
 Let $T$ be the tangent space to
$P_E$ at $o$ and let $T', T'' \subset T$ be the images of $du_E$ at $o'$, $o''$.
Since $du_E$ is injective at $o', o''$ and $u_E^{-1}(o) = \lbrace o', o''
\rbrace$, it follows that $T' \cup T''$ spans $T$ and that $T' \cap T''$ is the
tangent space to $\sing \ P_E$ at $o$. We have $\dim T' \cap T'' \geq 1$ because
$\sing \ P_E$ is a curve. On the other hand $P_{E,x} \cap P_{E,y} = o$ implies 
$\dim \ T \ \geq 2r-1$. Since $\dim \ T' = \dim \ T'' = r$, we deduce that $\dim
\ T' \cap T'' = 1$. Hence, as a scheme defined by the Jacobian ideal of
$P_E$, $\sing \ P_E$ is reduced. Finally the degree of $\sing \ P_E$ can be
obtained via double point formula, see \cite[9.3]{Fulton}, as follows:
\[
V \cdot V - c_{r-1}( {\emph N}_{V \vert {\P}^{2r-2}}) = 2 \ deg (\sing \ P_E)
\]
where $V \subset {\P}^{2r-2}$ is a general hyperplane section of $P_E$, corresponding to a global section $\sigma \in H^0({\mathcal O}_{\bold P_E}(1)  )\simeq H^0(\omega_C \otimes E)$, $c_{r-1}$ denotes the $(r-1)$- Chern class of the normal bundle ${\emph N}_{V \vert {\P}^{2r-2}}$. Note that  $V = P_{E'}$, with $E'$, vector bundle of rank $r-1$ defined by the section $\sigma$ as follows:
\[
0 \to {\mathcal O}_C \to E \otimes \omega_C \to E' \otimes \omega_C \to 0,
\]
and ${\mathcal O}_{\bold P_{E'}}(1)= {{\mathcal O}_{\bold P_E}(1)}_{\vert \P_{E'}}$.
Let $H = [{\mathcal O}_{\bold P_{E'}}(1)]$ and $f$ be the class of a fibre of  $\P_{E'}$, then 
 by computing  the total Chern class of the normal bundle, we find
\[
c_{r-1}( {\emph N}_{V \vert {\P}^{2r-2}})= r H^{r-1} + f H^{r-2}[4 r^2 -4r]= 7 r^2 - 4r.
\]
Finally,  we have $\deg \sing \ P_E = deg \ H_E = r^2 + 2r$. \par \noindent 
\end{proof}
 \par \medskip \noindent The genus of $D_E$ is $r^2+1$ and $H_E$ has degree
$r^2 + 2r$, hence $h^0(H_E) \geq 2r$
\begin{prop} \label{fourfive} For a general $E$ the line bundle $H_E$ is non special that is 
\[
h^0(H_E) = 2r. 
\]
\end{prop}
\begin{proof} By induction on $r$. Let $r = 1$ then $A := \omega_C \otimes E$
is a general line bundle of degree 3, $\P_E = C$ and
$u_E: \P_E \to \P^1$ is the triple covering defined by $A$.
Moreover $D_E$ is the family of divisors $x+y$ which are contained in a fibre
of $u_E$. It is easy to see that $D_E$ is a copy of $C$ and that $h_E = u_E$.
Then $H_E = A$ and hence $h^0(H_E) = 2$. \hfill \par \noindent Let $r \geq 2$
and let $[E_{r-1}] \in {\mathcal  U}_{r-1}$ and $E_1 \in \pic^1(C)$ be general points
satisfying the statement, then their corresponding curves $D_{r-1}$ and $D_1$
are smooth and transversal. Taking a general semistable extension
\begin{equation} 
0 \to E_{r-1} \to E \to E_1 \to 0 \label{succ}
\end{equation} 
we have $h^0(E \otimes \omega_C(-x-y)) \leq 1$ for any $x+y$,
(see \ref{twotwo} and its proof). Observe also that $D_E = D_1
\cup D_{r-1}$ and that $h_E$ is a morphism. The restrictions of $h_E$ to $D_1$
and $D_{r-1}$ can be described as follows: \hfill \par
\noindent (a) Let $F := \omega_C \otimes E$ and let $A := \omega_C \otimes
E_1$: tensoring \ref{succ} by $\omega_C$ and passing to the long exact
sequence, we obtain a surjective map  $H^0(F) \to H^0(A)$. Its dual is a linear
embedding $i: \P H^0(A)^* \to \P^{2r-1}$. On the other hand we
already know that $h_{E_1}$ is the triple cover of $\P H^0(A)^*$ defined
by $A$. It is easy to conclude that  $h_E|_{D_1} = i \cdot h_{E_1}$. In
particular
$h_E(D_1)$ is a line $\ell$ in $\P^{2r-1}$ which is triple for $h_E(D_E)$.
\hfill \par \noindent (b) Let $B := \omega_C \otimes E_{r-1}$: tensoring
(\ref{succ}) by $\omega_C$ and passing to the long exact sequence we get
an injection
$H^0(B) \to H^0(F)$. Its dual induces a projection $p: \P^{2r-1} \to
\P H^0(B)^*$ of centre $\ell$. It is again easy to conclude that $p \cdot
h_EE|_{D_{r-1}}$ $=$ $h_{E_{r-1}}$. \par \noindent It follows from the remarks
in (b) that
\[
 H_E \otimes{\mathcal  O}_{D_{r-1}} = H_{E_{r-1}} (a), 
\]
where $a := (h_E|_{D_{r-1}})^*\ell = D_1 \cdot D_{r-1}$. On the other hand (a)
implies that
\[
H_E \otimes {\mathcal  O}_{D_1} = H_{E_1} = A.
\]
Finally, tensoring by $H_E$ the Mayer-Vietoris exact sequence
\[
  0 \to {\mathcal  O}_{D_E} \to {\mathcal  O}_{D_{r-1}} \oplus {\mathcal 
O}_{D_1} \to {\mathcal  O}_a\to 0,
\]
we obtain
\begin{equation*} 0 \to H_E \to H_{E_{r-1}}(a) \oplus A \to {\mathcal  O}_a \otimes
H_E \to 0.
\end{equation*} 
By induction $h^1(H_{E_{r-1}}) = 0$, hence $h^1(H_{E_{r-1}}(a))
= 0$. Moreover $h^1(A) = 0$. Passing to the long exact sequence, the vanishing
of $H^1(H_E)$ follows if the restriction
\[
\rho: H^0(H_{E_{r-1}}(a)) \to {\mathcal  O}_a(a) \otimes H_{E_{r-1}}
\]
is surjective. Since $h^1(H_{E_{r-1}}) = 0$, this follows from the long
exact sequence of
\[
 0 \to H_{E_{r-1}} \to H_{E_{r-1}}(a) \to {\mathcal  O}_a(a) \otimes
H_{E_{r-1}}
\to 0.
\]
The vanishing of $h^1(H_E)$ extends by semicontinuity to a general point of
${\mathcal  U}_r$.
\end{proof}
\section{The Brill-Noether curve of $E$} In the following we will set for
simplicity: $D := D_E$. $D$ is an abstract curve endowed with an embedding $D
\subset C^{(2)}$. These data are in general not sufficient to reconstruct the
vector bundle $E$. As we will see the additional datum of $H_E$ makes possible
such a reconstruction. \medskip \par  \noindent The embedding $D \subset
C^{(2)}$  uniquely defines the family of divisors
\begin{equation*}
 b_x := C_x \cdot D
\end{equation*} 
where $x \in C$ and $C_x := \lbrace x+y\mid y \in C \rbrace$.
$b_x$ fits in the standard exact sequence 
\begin{equation*} 
0 \to H^0(\omega_C \otimes E(-x)) \otimes {\mathcal  O}_C \to
\omega_C \otimes E(-x) \to {\mathcal  O}_{b_x} \to 0
\end{equation*} 
and its degree is $2r$. The determinant of $E$ can be
reconstructed from the family $\lbrace b_x\mid   x \in C \rbrace$. Indeed let
$x+x'
\in | \omega_C |$, then the previous exact sequence implies 
\begin{equation*} 
\det \ E \cong {\mathcal  O}_C(b_x-rx').
\end{equation*} 
\rm \noindent Let $t+\Theta_E$ be the translate of $\Theta_E$ by
$t \in \mbox{\rm Div}^0(C)$:  $D_{E(-t)} = a^*(t+\Theta_E)$. Thus, up to
replacing $E$ by $E(-t)$, $D$ is transversal to $C_x$ and $b_x$ is smooth for a
general $x$. Mainly we will consider $b_x$ as a divisor on $D$. Let $d \in D$,
it is clear that: $d \in \mbox{\rm Supp} \ b_x$ $\Longleftrightarrow$ $d = x+y$
$\Longleftrightarrow$ $h_E(d) \in P_{E,x}$. This implies that
\begin{equation*} 
b_x \ = \ {h_E}^* P_{E,x}
\end{equation*} 
for each $x \in C$. The line bundles $H_E(-b_x)$ have degree
$r^2$. Since $D$ has genus $r^2+1$ they define a family of points in the theta
divisor of $\pic^{r^2}(D)$. We can say more:
\begin{prop}\label{fiveone}  Let $E$ be general then $h^0(H_E(-b_x)) = r$, for
each $x\in C$.
\end{prop}
\begin{proof} We know from Prop.~\ref{fourfive} that $h^0(H_E) = 2r$, we also know that
$b_x$ $=$ ${h_E}^*(P_{E,x})$. Since the space $P_{E,x}$ has dimension $r-1$, it
follows $h^0(H_E(-b_x)) \geq r$. Moreover the equality holds if the set
$h_E(\mbox{\rm Supp} \ b_x)$ spans $P_{E,x}$. We prove this by induction on $r$. Let $r =
1$ then $P_{E,x}$ is a point: since $h_E$ is a morphism $h_E(\mbox{\rm Supp} \ b_x) =
P_{E,x}$. Let $r \geq 2$, as in the proof of \ref{fourfive} we consider a general extension
\begin{equation} \label{nineteen}
0 \to E_{r-1} \to E \to E_1 \to 0
\end{equation} 
with $[E_{r-1}] \in {\mathcal  U}_{r-1}$ and $E_1 \in {\pic}^1(C)$
general. We fix the same assumptions and notations of the proof of
\ref{fourfive} which is similar. In particular the curves $D_{E_{r-1}}$ and
$D_{E_1}$ are smooth and transversal, moreover the exact sequence
\[
 0 \to B \to F \to A \to 0
\]
just denotes the above exact sequence (\ref{nineteen}) tensored by
$\omega_C$.  Such a sequence induces a linear embedding $i: \P H^0(A)^* \to
\P H^0(F)^*$. The image of $i$ is the line $\ell$ considered in \ref{fourfive}
and it holds the equality proved there: $h_E|_{D_{E_1}}$ $=$
$i \cdot u_{E_1}$. Then it turns out that
\[
 \ell \cap P_{E,x} = h_E(C_x \cap D_{E_1}) = \ \text {one point $o_x$}
\]
for each $x \in C$. On the other hand let $p: P_E \to \P H^0(B)^*$ be
the projection of centre $\ell$, then the latter exact sequence implies that
$p(P_{E,x}) = P_{E_{r-1},x}$. Moreover we also know from the proof of
\ref{fourfive} that
$h_{E_{r-1}} = p \cdot h_E$. By induction 
$P_{E_{r-1},x}$ is spanned by 
\[
 h_{E_{r-1}}(C_x \cap D_{E_{r-1}}) = p(h_E(C_x \cap D_{E_{r-1}})). 
\]
 Hence
the linear span of $h_E(C_x \cap D_{E_{r-1}})$ is a space
$L \subset P_{E,x}$ of dimension $\geq r-2$ and such that $p(L) =
P_{E_{r-1},x}$. If $o_x \in P_{E,x} - L$ then $P_{E,x}$ is spanned by $h_E(C_x
\cap D)$. If
$o_x \in L$ then $\dim L = r-1$ and $L = P_{E,x}$. In both cases the statement
follows.
\end{proof} \medskip \par \noindent For each $\ell \in \pic^1(C)$ we consider the
curve 
\[
 B_\ell := \lbrace x+y \ \in \ | \ell(z) |, \ z \in C \rbrace.
\]
$B_\ell$ is biregular to $C$ unless $\ell = {\mathcal  O}_C(x)$ for some $x \in C$. In the
latter case $B_\ell$ is  $C_x \cup \ | \omega_C |$. We define
\begin{equation*} 
b_\ell := D \cdot B_\ell.
\end{equation*} 
Note that $b_\ell = b_x$ if $B_\ell = C_x \cup |\omega_C|$.
The reason is that we are assuming $E$ general, then $h^0(E) = 0$ and hence
$D \ \cap | \omega_C | = \emptyset$. 
\lem \label{new} The morphism $ b: \pic^1(C) \to \pic^{2r}(D)$ sending $\ell$ to $b_\ell$ is an embedding. We will denote its image as $J_D$:
\begin{equation*} 
J_D := \lbrace {\mathcal  O}_D(b_\ell)\mid  \ell \in \pic^1(C)
\rbrace,
\end{equation*}  
in particular $J_D$ contains the canonical theta divisor
\begin{equation} \label{twentytwo} C_D := \lbrace {\mathcal  O}_{D}(b_x)\mid  \
x \in C
\rbrace.
\end{equation} 
\endrem
\par
\noindent
\endlem \rm
\begin{proof} 
Up to shifting the degrees, $b$ is a morphism between the complex tori $\pic^0(C)$
and $\pic^0(D)$. Hence it is an isogeny up to translations, so  $b$ is an embedding   if it  is 
injective. Let $\ell_1, \ell_2 \in \pic^1(C)$ and set \ \ \ ${\mathcal  L} \colon = {\mathcal  O}_D(b_{\ell_1}- b_{\ell_2})$.  ${\mathcal  L}$ is defined by the standard exact sequence:
\[
0 \to {\mathcal  O}_{C^{(2)}}(- D + B_{\ell_1}- B_{\ell_2}) \to  {\mathcal  O}_{C^{(2)}}( B_{\ell_1}- B_{\ell_2}) \to
{\mathcal  L} \to 0. 
\]
$D + B_{\ell_2}- B_{\ell_1}$ is the pull-back  by the Abel map $a \colon C^{(2)} \to \pic^0(C)$ of a divisor
homologous to $r\Theta$, where $\Theta$ is a theta divisor in $\pic^0(C)$. Since $r\Theta$ is ample, it follows:
$h^0(- D + B_{\ell_1}- B_{\ell_2}) = h^1( - D + B_{\ell_1}- B_{\ell_2}) = 0$.  
So the associated long exact sequence gives:
\[
h^0({\mathcal  L} ) = h^0({\mathcal  O}_{C^{(2)}}( B_{\ell_1}- B_{\ell_2})).
\]
Moreover, it's easy to see that if ${\ell_1}\not= {\ell_2}$, then  $h^0({\mathcal  O}_{C^{(2)}}( B_{\ell_1}- B_{\ell_2})) = 0$
 hence $b_{\ell_1}$ and $b_{\ell_2}$ are not linearly equivalent.
\end{proof}
As an immediate consequence of the lemma, the following map 
\[
b_0 \colon \pic^0(C) \to \pic^0(D)
\]
  sending ${\mathcal  O}_C(x-y)  \to {\mathcal  O}_D(b_x-b_y)$ is an embedding too.
\hfill\par
\hfill\par
As we already pointed out $C_D$ is not sufficient to reconstruct
$E$, the crucial curve for doing this can be now defined:
\defi The Brill-Noether curve of $E$ is the curve
\begin{equation*} C_E := \lbrace H_E(-b_x), \ x \in C \rbrace.
\end{equation*}
\enddefi \rm \noindent
$C_E$ is a copy of $C$ embedded in $\pic^{r^2}(D)$. Since $h^0(H_E(-b_x)) = r$,
each point of $C_E$ is a point of multiplicity $r$ for the theta divisor
\[
\Theta_D := \lbrace L \in \pic^{r^2}(D) \mid  h^0(L) \geq 1 \rbrace.
\]
In particular $C_E$ is contained in the Brill-Noether locus
\[
W^{r-1}_{r^2}(D) := \lbrace L \in \pic^{r^2}(D) \mid h^0(L) \geq r \rbrace.
\]
The  Brill-Noether number $\rho(r-1,r^2,r^2+1)$ yields the expected dimension
of $W^{r-1}_{r^2}(D)$. We have 
$ \rho(r-1,r^2,r^2+1) = 1$ for each $r$, so we expect that $C_E$ is an
irreducible component of $W^{r-1}_{r^2}(D)$, see \cite[ch. V]{Arba}. Of course
$D$ is not a general curve of genus $r^2+1$, so the latter property is not a
priori granted. 
\hfill\par
\hfill\par
\begin{rem} \label{fivethree} $E$ is uniquely reconstructed from the pair
$(C_D,H_E)$ as follows:
  \par \noindent
\proof Consider the correspondence
\[
 B = \lbrace \ (x,y+z) \in C \times D \mid x \in \lbrace y, z \rbrace \
\rbrace,
\] 
with  $B \cdot (x \times D )= b_x$. Let $p_1: C \times D \to C$
and $p_2: C \times D \to D$ be the projection maps, then we apply the functor
$p_{1*}$ to the exact sequence
\[
 0 \to p_2^*H_E(-B) \to p_2^*H_E \to p_2^*H_E \otimes {\mathcal  O}_B \to 0.
\]
This yields the exact sequence
\[
 0 \to \overline F^* \to H^0(H_E) \otimes {\mathcal  O}_C \to p_{1_*}
{\mathcal  O}_B
\otimes p_2^*H_E \to R^1p_{1_*}p_2^*H_E(-B) \to 0,
\]
 where 
\[\overline F^* := p_{1*} p_2^*H_E(-B).
\]
 Let $F = \omega_C \otimes
E$, we have the natural identities
\[
 {\overline {F_x}^*} = H^0(H_E(-b_x)) = H^0(F(-x)).
\]
The left one is immediate. Let $\mathcal I$ be the Ideal of $P_{E,x}$, then
we  have
$H^0(H_E(-b_x)) = H^0({\mathcal  I}(1)) \ $ by prop. 5.1. Hence the right equality follows from the
identity $H^0({\mathcal  I}(1))$ $=$ $H^0(F(-x))$. The above identities, together
with $H^0(H_E) = H^0(F)$, imply that
\[
 F = H^0(H_E) \otimes {\mathcal  O}_C / \overline F^*.\quad\quad  \qed
\]
\end{rem} \noindent As an immediate consequence of the above
construction we have:
\begin{prop} \label{fivefour} Let $[E_1], [E_2]$ be general points of 
${\mathcal  U}_r$. Assume that:
${\theta}_r([E_1]) = {\theta}_r([E_2]) = D$ and   $H_{E_1}= H_{E_2}$. Then
$[E_1] = [E_2]$.
\end{prop}
\rem \label{fivefive} \rm The previous construction also defines the vector
bundles
\begin{equation*}
\overline E \colon = \overline F \otimes {{\omega}_C^{-1}} \quad , \quad \widetilde E
\colon =  {{\omega}_C^{-1}} \otimes (R^1{p_1}_* p_2^*H_E(-B)).
\end{equation*} 
We already know from \ref{threetwo} that the assignment $[E] \to
[\overline E]$ defines a birational involution $j: {\mathcal  U}_r \to {\mathcal  U}_r$.
Notice also that
$\widetilde E$ is semistable for $E$ general: to prove this it suffices to produce
one semistable $E_o$ such that $\widetilde E_o$ is semistable. The existence of 
$E_o$ follows by induction on $r$: this is obvious for $r = 1$. Let $r \geq 2$
and let $E_o$ be defined by a semistable extension $e \in \mbox{\rm
Ext}^1(E_1,E_{r-1})$ where $[{\mathcal  U}]_{r-1} \in {\mathcal  U}_{r-1}$ and $E_1
\in {\mathcal  U}_1$. It is easy to show that $\widetilde E_o$ is defined by some $\widetilde e
\in \mbox{\rm \mbox{\rm \mbox{\rm \mbox{\rm Ext}}}}^1(\widetilde E_1,
\widetilde E_{r-1})$: we leave the details to the reader. Hence $\widetilde E_o$ is
semistable. Due to this property we can define a rational map
\[
\kappa: {\mathcal  U}_r \to {\mathcal  U}_r
\]
sending $[E]$ to $[\widetilde E]$. In addition we have:
\endrem
\begin{prop} \label{fivesix} $\kappa: {\mathcal  U}_r \to {\mathcal  U}_r$ is
birational, in particular its inverse is $j \cdot \kappa \cdot j$.
\end{prop}
\begin{proof} Let $T := {\mathcal  O}_D(b_x+b_{i(x)})$, where $i:C \to C$ is the
hyperelliptic involution. $T$ does not depend on $x$ because the family of
divisors $ \lbrace b_x+b_{i(x)}, \ x+i(x) \in \ | \omega_C | \rbrace$ is
rational. Then we define the line bundle of degree
$r^2 + 2r$
\[
\widetilde H_E := \omega_D \otimes T \otimes H_E^{-1}.
\]
First, we note that $h^1( \widetilde H_E) = 0$ for $E$ general. Indeed
$\omega_D \otimes \widetilde H_E^{-1}$ is $H_E(-b_x - b_{i(x)})$ and hence
$h^1(\widetilde H_E) = h^0(H_E(-b_x-b_{i(x)})$ by Serre duality. Since
$h^0(H_E(-b_x-b_{i(x)})$ $=$ $h^0(\omega_C \otimes E(-x-i(x))$ $=$ $h^0(E)$, it
follows $h^1(\widetilde H_E) = 0$ for each $[E] \in {\mathcal  U}_r - \Theta_r$. Secondly
we note that, with the previous notations, Serre duality yields a natural
identification
\[
 R^1p_{1_*}p_2^*H_E(-B)_x = H^0(\widetilde H_E(-b_{i(x)}))^*, \ \forall x \
\in \ C.
\]
It is then easy to deduce that
\[
\omega_C \otimes \widetilde E = R^1p_{1_*}p_2^*H_E(-B) \cong p_{1*}p_2^* \widetilde
H_E(-B)^*.
\]
 Starting from $\widetilde H_E$ it is clear that one obtains $H_E$ and with
the same construction $\overline E = \omega_C ^{-1} \otimes p_{1*}p_2^*H_E$. 
 Notice also that $\widetilde H_E$ is the line bundle $H_{\widetilde
{\overline E}}$ defined by the vector bundle $\overline {\widetilde E}$. This implies
that $\kappa^{-1}$
$=$ $j \cdot \kappa \cdot j$ and hence that $\kappa$ is birational.
\end{proof}
\section {The fibres of the theta map}  We want to see that $E$ is also
uniquely reconstructed from the pair $(D,C_E)$. For this we consider more in
general any smooth curve $D
\subset C^{(2)}$ such that $a_*D \in {\mathcal  T}_r$:
\defi \label{sixone} A Brill-Noether curve of $D$ is a copy 
\[
 C' \subset \pic^{r^2}(D)
\]
 of
$C$ satisfying the following  property: there exists $H \in \pic^{r^2 + 2r}(D)$
such that
\[
 C' = \lbrace H(-b_y), y \in C \rbrace,
\]
moreover $H$ is non special and $h^0(H(-b_x)) = r$ for every point $x \in
C$. The set of the Brill-Noether curves of $D$ will be denoted by $S_D.$
\enddefi \rm \par \noindent
 Let $D = \theta_r([E])$ then $C_E$ is a Brill-Noether curve of $D$. Let $r =
1$ then $D = C$ and the canonical theta divisor of
$\pic^1(C)$ is the unique Brill-Noether curve of $D$. \noindent
\lem \label{sixone} Let $H$ be as in the previous definition then $H$ is unique.
\endrem
\par
\noindent
\endlem \rm
\begin{proof} Assume that $C' = \lbrace H'(-b_x), x \in C \rbrace$ for a second
$H'$. Then there exists an automorphism $u: C \to C$ which is so defined $u(x)
= y$ iff $H'(-b_x)$ $=$ $H(-b_y)$. Let $\gamma: C \times C \to \pic^0(D)$ be the
map sending $(x,y)$ to $H'
\otimes H^{-1} (b_x-b_y)$: we will show  that the image of $\gamma$ is a copy  of $\pic^0(C)$ and that 
$\gamma: C \times C \to \pic^0(C)$ is the  difference map. 
To see this recall
that $C_D = \lbrace b_x, \ x \in C \rbrace$ is the theta divisor of $J_D$, see
(\ref{twentytwo}). The map $\tilde{\gamma} \colon C \times C \to \pic^0(D)$ sending
$ (x,y)$ $\to$ ${\mathcal O}_D(b_x-b_y)$ factors through the isomorphism $t \colon C \times C \to C_D \times C_D$, sending
$(x,y) \to (b_x, b_y)$, and  the difference map.  Moreover, we have the following commutative diagram
\[
\begin{array}{c}
{C \times C} \\ {\downarrow} \\  \pic^0(C) \end{array}
\begin{array}{c}
\stackrel{t}{\rightarrow} \\   \\  \stackrel{b_0}{\rightarrow} \end{array}
\begin{array}{c}
{C_D \times C_D} \\ {\downarrow} \\  \pic^0(D) \end{array}
\]
where the vertical arrows are difference maps and $b_0 \colon \pic^0(C) \to \pic^0(D)$
sending ${\mathcal O}_C(x-y) \to {\mathcal O}_D(b_x-b_y)$ is an embedding, see \ref{new}.
This implies  that $\tilde{\gamma}$ is a difference map and  $\gamma$ too.
  The graph of $u$ is obviously contracted by $\gamma$, on the other
hand the only curve contracted by the difference map is the diagonal of $C
\times C$. Then
$u$ is the identity and $H(-b_x) = H'(-b_x)$ for each $x \in C$. Hence $H = H'$.
\end{proof} \noindent
\begin{prop} \label{sixthree} Let $[E_1], [E_2]$ be general points in $U_r$,
assume that
$C_{E_1}= C_{E_2}$ and that 
${\theta}_r([E_1]) = {\theta}_r([E_2]) = D.$ Then $[E_1] = [E_2]$.
\end{prop}
\begin{proof} By the previous lemma  $H_{E_1}= H_{E_2}$ and this implies $[E_1]
= [E_2]$, by  \ref{fivefour}.
\end{proof}
\teo \label{sixfour} The theta map $\theta_r: {\mathcal  U}_r \to {\mathcal 
T}_r$ is generically finite.
\endteo
\begin {proof} It suffices to show that $\theta_r|_U: U \to \theta_r(U)$ is
generically finite for a suitable dense open set $U$. Hence we can assume that
$D \in \theta_r(U)$ is a smooth curve and that the points of
${(\theta_r|_U)}^{-1}(D) = \theta_r^{-1}(D) \cap U$ are sufficiently general in
${\mathcal  U}_r$. Let
\[
 i: (\theta_r|_U)^{-1}(D) \to S_D
\]
be the map sending $[E]$ to $C_E$. By \ref{sixthree} $E$ is uniquely
reconstructed from
$(D,C_E)$, hence it follows that $i$ is injective. On the other hand recall
that $C_E$ is contained in the Brill-Noether locus $W^{r-1}_{r^2}(D)$. Since
the Brill-Noether number
$\rho(r-1,r^2,r^2+1)$ is one, each irreducible component of $W^{r-1}_{r^2}(D)$
has dimension $\geq 1$. This implies  that $\theta_r|_U$ is finite if $C_E$ is
an irreducible component of $W^{r-1}_{r^2}(D)$. This property is proved in the
next theorem. Hence the statement follows. 
\end{proof}
\begin{lem}  \label{sixfive} Let $D = \theta_r([E])$ for a general $[E] \in
{\mathcal  U}_r$ and let
$a = D \cdot D_1$ for a general $D_1 \in {\mathcal  T}_1$, then the line bundle
$H_E(a-b_x)$ is non special.
\end{lem}
\begin{proof} Let  $D_1 = \theta_1([E_1]) \subset C^{(2)}$ with $E_1 = {\mathcal 
O_C}(x)$, then we have: $D_1 = x \times C \ \cup | \omega_C | \subset C^{(2)}$.
Note that  
$a =D
\cdot D_1= b_x$ if
$E$ is general. Hence
$H_E(a-b_x) = H_E$ is non special. By semicontinuity, the same is true for a
general $D_1$.
\end{proof}
\begin{teo} For a general $[E] \in {\mathcal  U}_r$ the Brill-Noether curve $C_E$ is
an irreducible component of $W^{r-1}_{r^2}(D)$, $D = {\theta}_r([E])$.
\end{teo}
\begin{proof} Let $H := H_E$, it is sufficient to show the injectivity of the
Petri map
\begin{equation*}
\mu: H^0(H(-b_x)) \otimes H^0( \omega_D \otimes H^{-1}(b_x)) \to H^0(\omega_D)
\end{equation*} 
for a general $x \in C$. This implies that the tangent space to
$W^{r-1}_{r^2}(D)$ at its point $H(-b_x)$ is 1-dimensional,  see
\cite[Ch.~V]{Arba}. We proceed by induction on $r$. Let $r = 1$, then $D =
C$ and $C_E =
\lbrace {\mathcal  O}_C(x), x \in C \rbrace$. Hence the injectivity of
$\mu$ is immediate. For $r \geq 2$ we borrow once more the notations and the
method from the proof of proposition \ref{fourfive}. So we specialize
$E$ to the semistable vector bundle defined by the exact sequence
\[
 0 \to E_{r-1} \to E \to E_1 \to 0.
\]
Then $D$ is the transversal union of the curves $D_{r-1} =
\theta_{r-1}([E_{r-1}])$ and
$D_1 = \theta_1([E_1])$, $h_E$ is a morphism and $H$ is the line bundle $h_E^*
O_{{\bf P}^{2r-1}}(1)$. Let $a = D_{r-1} \cdot D_1$: from \ref{fourfive}  we have $D_1 = C$ and $H_{E_1}= \omega_C \otimes E_1$ and  moreover
\[
 H \otimes {\mathcal  O}_{D_{r-1}} = H_{E_{r-1}}(a) \ \text {and}  \ H \otimes
{\mathcal  O}_{D_1} = H_{E_1}.
\]
 Since $x$ is general we can assume $\mbox{\rm Supp} \ b_x \cap \sing \ D =
\emptyset$ so that ${\mathcal  O}_D(b_x)$ is a line bundle. Let $\mathcal I$ be
the ideal sheaf of
$D_1$ in $D$: at first we show that $\mu  \vert  I \otimes W$ is injective,
where
\[
 I := H^0({\mathcal  I} \otimes H(-b_x)) \ \text {and} \ W := H^0(\omega_D
\otimes H^{-1}(b_x)).
\]
Since ${\mathcal  I} \otimes {\mathcal  O}_{D_{r-1}}$ $=$ ${\mathcal 
O}_{D_{r-1}}(-a)$ and
$\omega_D \otimes {\mathcal  O}_{D_{r-1}}$ $=$ $\omega_{D_{r-1}}(a)$, we have the
restriction maps
\[
\rho_I: I \to H^0(H_{E_{r-1}}(-b_{x,r-1})) \ \text {and} \ \rho_W: W \to
H^0(\omega_{D_{r-1}} \otimes H_{E_{r-1}}^{-1}(b_{x,r-1}))
\]
with $b_{x,r-1} = b_x \cdot D_{r-1}$. 
\hfill\par {\bf Claim }: $\rho_I$ is an isomorphism and  $\rho_W$ is
surjective. 
\hfill\par Let's assume the claim, then $\rho := \rho_I \otimes \rho_W$ is
surjective and defines the exact sequence
\[
\hspace{-1em} 0 \to \ker \ \rho \to I \otimes W \to H^0(H_{E_{r-1}}(-b_{x,r-1}))
\otimes H^0(\omega_{D_{r-1}} \otimes H_{E_{r-1}}^{-1}(b_{x,r-1})) \to 0.
\]
In particular it follows $\dim \ \ker \ \rho = r-1$. By induction on $r$ the
Petri map on the tensor product at the right side is injective. Therefore $\mu
\vert I \otimes W$ is injective iff $\mu  \vert \ker \ \rho $ is injective. But
our claim implies $\dim \ \ker \ \rho_W = 1$ and $\ker \ \rho$ $=$ $I \
\otimes \langle w\rangle $, where $w$ generates $\ker \ \rho_W$. Hence $\mu 
\vert  \ker
\
\rho$  is injective as well as $\mu \vert  I \otimes W$. \hfill \par \noindent 
Let $V := H^0(H(-b_x))$, now we consider the exact sequence
\[
 0 \to I \otimes W \to V \otimes W \to (V/I) \otimes W \to 0.
\]
The map $\mu$ induces a multiplication 
\[
\nu: (V/I) \otimes W \to
H^0(\omega_D)/\mu(I \otimes W).
\]
 The injectivity of 
$\mu \vert  I \otimes W$ implies that $\mu$ is injective iff $\nu$ is
injective. On the other hand, $\rho_I$ is an isomorphism, hence  
$\dim \ I = r-1$ and $\dim \ V/I = 1$. Let $v \in V-I$ then: $\nu$ is injective
iff $vW \cap \mu(I \otimes W) = (0)$ iff no $w \in W-(0)$ vanishes on $D_1$.
This is equivalent to the injectivity of the restriction map
\[
 u: H^0(\omega_D \otimes H^{-1}(b_x)) \to H^0(\omega_C(a-x));
\] 
in fact $W = H^0(\omega_D \otimes H^{-1}(b_x))$ and $\omega_D
\otimes H^{-1}(b_x) \otimes {\mathcal  O}_{D_1} = \omega_C(a-x)$. To prove that
$u$ is injective consider the Mayer-Vietoris long exact sequence
\[
 0 \to W \to H^0(\omega_{D_{r-1}} \otimes H_{E_{r-1}}^{-1}(b_{x,r-1})) \oplus
H^0(\omega_C(a-x)) \to {\mathcal  O}_a \to \dots
\]
  The left non zero arrow followed by the projection onto
$H^0(\omega_C(a-x))$ is exactly $u$. This implies that $\ker \ u$, via the
restriction map, injects  in
$H^0(\omega_{D_{r-1}} \otimes H_{E_{r-1}}^{-1}(b_{x,r-1}-a))$.
 So $u$ is injective if the latter space is zero that is if
$H^1(H_{E_{r-1}}(a-b_{x,r-1})) = 0$: this has been shown in lemma~\ref{sixfive}.
Hence
$\mu$ is injective. Then, by semicontinuity, the same property is true for a
general $[E] \in {\mathcal  U}_r$ and the statement follows.
\hfill\par \noindent To complete the proof we show the above  claim. \hfill\par
\noindent - Let $h: D \to \P^{2r-1}$ be the map defined by $H$. As in
\ref{fourthree}  $h(D_1)$ is a line $\ell$ and $P_{E,x} \cap \ell$ is a point.
Moreover
$P_{E,x}$ is spanned by
$h(b_x)$. Hence we have $I = H^0({\mathcal  J}(1))$ and $\dim \ I = r-1$, ${\mathcal  J}$
being the ideal of
$P_{E,x} \cup \ell$. In particular $\rho_\ell$ is the pull-back
$(h|_{D_{r-1}})^*$ restricted to a space of linear forms vanishing on $h(D_1)$.
Since $h(D)$ is non degenerate $\rho_\ell$ is injective. Then, for dimension
reasons,
$\rho_\ell$ is an isomorphism.
\hfill\par \noindent - As in \ref{fourfive} the projection $p: P_{E,x} \to
P_{E_{r-1},x}$ from $P_{E,x} \cap \ell$ is surjective. Equivalently the
restriction
$H^0(E\otimes \omega_C(-x)) \to H^0(E_{r-1}\otimes \omega_C(-x))$ is
surjective. So this property holds for general $E,E_{r-1}$. Let $\widetilde E_r$,
${\widetilde E}_{r-1}$ be defined from $E_r$, $E_{r-1}$ as in
Remark~\ref{fivefive}. By
\ref{fivesix} they are general. Hence the restriction $H^0(\widetilde
E\otimes
\omega_C(-x))
\to H^0({\widetilde {E}}_{r-1}\otimes
\omega_C(-x))$ is surjective: this map is just $\rho_W$. 
\end{proof}
\hfill\par\rm
\noindent We can summarize as follows our partial geometric description of the
theta map:
\begin{teo} \label{sixseven} Let $D \in {\mathcal  T}_r$ be general and smooth,
then there exists a natural injective map $i_D$ between the fibre of $\theta_r$
at
$D$ and the set of the Brill-Noether curves of $D$. Namely the map
\[
i_D: {\theta_r}^{-1}(D) \to S_D
\]
associates to $[E] \in \theta_r^{-1}(D)$ its Brill-Noether curve $C_E \in
S_D$. 
\end{teo} 
\begin{proof} Since $\theta_r$ is generically finite, each point $[E]$ $\in$
$\theta^{-1}(D)$ is sufficiently general in ${\mathcal  U}_r$. The injectivity then
follows from corollary \ref{sixfour}.
\end{proof}
\hfill\par \begin{rem} \rm Each Brill-Noether curve $C \in S_D$ uniquely
defines a vector bundle $E_C$ of rank $r$ and degree $r$:  to construct $E_C$
it suffices to take the line bundle $H$ appearing in the definition \ref{sixone}
of Brill-Noether curve. Applying to the pair $(C_D,H)$ the reconstruction
produced in remark \ref{fivethree} we finally obtain such a vector bundle
$E_C$. If $E_C$ is semistable it turns  out  that $\theta_r([E_C]) = D$ and
that $i_D([E_C]) = C$. In particular $i_D$ is bijective if each $C \in S_D$
defines a semistable
$E_C$. This property seems very plausible for a general $D$, however we do not
have a rigorous proof of it.
\end{rem}

\end{document}